# DYNAMIC INSPECTION PLANNING FOR SYSTEMS WITH INDIVIDUALLY REPAIRABLE COMPONENTS

NOOSHIN YOUSEFI, DAVID W. COIT

*Industrial and System Engineering Department, Rutgers Univeristy*
*Piscataway, NJ, 08854, USA*
*Email: no.yousefi@rutgers.edu, Coit@soe.rutgers.edu*

In this paper a multi-component system is studied where each component can be repaired within the system. We consider that each component is subject to two dependent competing failure processes due to degradation and random shocks and each shock increases the total degradation of each component. For systems with individual repairable components, each failed component can be replaced at any inspection time and other components continue functioning, so the age of each component at any inspection time is different from other components. Different initial ages have an effect on the optimal time that the whole system should be inspected, therefore the optimal inspection time should be calculated dynamically considering the age of each component at the beginning of the interval. In this paper, the system reliability and cost rate function is calculated for a system with individually repairable components considering random age for each component. Two numerical examples demonstrate the proposed reliability and maintenance model.

## 1. Introduction

There are many systems that fail due to competing failure processes such as internal degradation and random shocks simultaneously. Failure due to continuous degradation is known as soft failure, and failure due to instantaneous stress caused by random shocks is called hard failure. It is assumed that each arriving shock may have some damages as additional abrupt change on the total degradation.

In the real word, the downtime cost is often much higher than the cost of repairing/replacing the system or components. Therefore, implementing a suitable preventive maintenance policy would be a cost saving plan. Preventive maintenance functions can restore the system to the state before failure by repairing or replacing the system or components to avoid failure.

In this paper, we study systems with individually repairable components, where each component can be repaired whenever it is failed. Previous research considered a system packaged and sealed together, and the whole system should be repaired/replaced upon failure, however for multi-component system with individually repairable components, it is not beneficial or cost effective to replace the whole system. So, developing a model for the reliability and maintenance for such systems with individually repairable components is a unique challenge to analyze. In this paper, we develop a reliability and maintenance model for a multi-repairable components system when each component is subject to two competing failure processes of hard and soft failure.

## 2. Background

There have been several studies on reliability and maintenance of multi-component systems subject to degradation and random shocks. Song et al. [1] developed a reliability model for a multi-component system experiencing multiple competing processes. Wang and Pham [2] studied the relationship of random shocks and a degradation process on the system reliability. Rafiee et al [3] developed a reliability model for a system experiencing degradation and random shocks with changing of degradation rate according to random shock arrivals.

There are different maintenance policies and techniques for single or multi-component systems. There have been studies on maintenance of single or multi-component systems subject to degradation and random shocks [4,5]. Abdul-Malak and Kharoufeh [6] used a Markov decision process model to develop a maintenance model for a multi-component system. Song et al. [7] developed a new condition-based maintenance model for a system of multi-components subject to dependent competing failure processes.

In this paper, we use a replace on failure preventive maintenance model for multi-repairable component system exposed to dependent competing failure processes.

## 3. System reliability analysis

In this paper, we consider two competing failure processes, soft failure and hard failure. Soft failure caused by continuous degradation and additional abrupt degradation damage by a shock process, and hard failure caused by instantaneous stress on system from each shock. Soft failure happens when the total degradation of component $i$ is greater than its own predefined soft failure threshold ($H_i$). When any shock magnitude for component $i$ is greater than a predefined hard threshold ($D_i$), component $i$ experiences hard failure. In this study, the Poisson process is used to model the shock arrival process. So, the probability of having $m$ shocks by time $t$ is given as follow:

$$P(N(t) = m) = \frac{(\lambda_0 t)^m e^{-\lambda_0 t}}{m!} \quad (1)$$

It is assumed that each shock magnitude is a *i.i.d* random variable following a normal distribution. $W_{ij} \sim Normal(\mu_{Wi}, \sigma^2_{Wi})$, where $W_{ij}$ is the $j^{th}$ shock magnitude of component $i$. Therefore, the probability of having no hard failure for component $i$ by time $t$ can be calculated by equation (2)

$$P_{NH,i}(t) = P(W_{ij} < D_i) = F_{W_i}(D_i) = F(\frac{D_i - m_{W_i}}{S_{W_i}}) \quad (2)$$

To calculate the probability of having no soft failure, we need to model the degradation process of each component. For a multi-component system each component degrades in the form of cumulative damage, and a stochastic process can be used to model the degradation process. In this paper, gamma process is used, which is a suitable stochastic process for monotonic increasing degradation path; So, $X(t)-X(s) \sim gamma(\alpha(t)- \alpha(s), \beta)$, where $X(t)$ is the degradation level at time $t$, $\alpha(t)=\alpha.t$ is gamma shape parameter which is linear in $t$, and $\beta$ is scale parameter.

Moreover, whenever each shock comes to the system, it has an additional abrupt damage on the total degradation of each component. It is assumed that the shock damage on each component is a *i.i.d* random variable which follows a normal distribution, $Y_{ij} \sim Normal(\mu_{Yi}, \sigma^2_{Yi})$ where $Y_{ij}$ is the $j^{th}$ shock damage on the $i^{th}$ component. $S_i(t)$ is the summation of all shock damages on component $i$ by time $t$. So, the probability that component $i$ has not experienced soft failure by time $t$ is given by equation (3)

$$P_{NS}(t) = P(X_{S_i}(t) < H_i) = \sum_{m=0}^{\infty} P(X_{S_i}(t) + S_i(t) < H_i) \times \frac{(\lambda_0 t)^m e^{-\lambda_0 t}}{m!} \quad (3)$$

In a multi-component system with individually repairable components, whenever each component experiences soft or hard failure, it is considered as failed and it should be replaced with a new one, but other components continue functioning until they fail. Since we replace the failed component instead of the whole system, the age of components at each inspection time is different from others. In this maintenance policy, the inspection interval should be found dynamically based on the initial age of all the components. If all components are new, the optimal inspection should be larger than the case that all components are prone to failure. For a system with multi-components that degrading differently, a preventive maintenance model should be found considering the age of all the components at the beginning of each interval.

Figure 1 shows the maintenance model for a multi-component system where each component can be replaced individually. At the beginning of each inspection interval, the ages of components are different from each other and subsequently the length of next inspection interval would be different based on the initial age of all the components. To calculate the system reliability for a future inspection interval, random values $u_i$ are assumed as the initial age of each component $i$ at the beginning of interval.

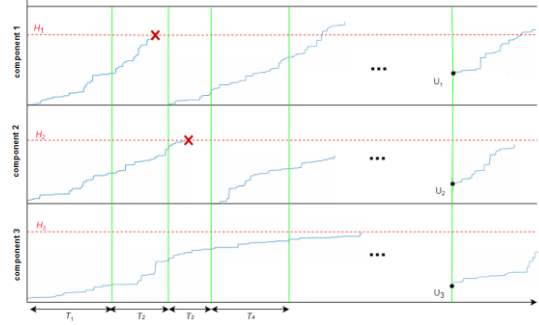

Figure 1. Dynamic inspection planning considering random initial age for each component

Equation (4) shows the probability of having no failure for component $i$ in a series system from the beginning of an interval up to time $t$.

$$R_i(t;u_i) = \sum_{m=0}^{\infty} \left[ P(W_{ij} < D_i)^m P\left(X_i(t) + S_i(t) + u_i < H_i\right) \mid N(t) = m \right] \frac{(\lambda_0 t)^m e^{(-\lambda_0 t)}}{m!} \quad (4)$$

$$= \sum_{m=0}^{\infty} \left[ P(W_{ij} < D_i)^m \times \int_0^{H_i} P\left(X_i(t) + y + u_i < H_i\right) f_{Y_i}^{<m>}(y) dy \right] \frac{(\lambda_0 t)^m e^{(-\lambda_0 t)}}{m!}$$

For a series configuration, the system fails if any component fails; hence, the system reliability can be calculated using equation (5), and for parallel configuration, the system fails if all the components are failed, so the system reliability should be calculated using equation (6)

$$R_S(t;\boldsymbol{u}) = \sum_{m=0}^{\infty} \prod_{i=1}^{n} \left[ P(W_{ij} < D_i)^m \times \int_0^{H_i} P\left(X_i(t) + y + u_i < H_i\right) f_{Y_i}^{<m>}(y) dy \right] \frac{(\lambda_0 t)^m e^{(-\lambda_0 t)}}{m!} \quad (5)$$

$$R_P(t;\boldsymbol{u}) = 1 - \sum_{m=0}^{\infty} \prod_{i=1}^{n} \left[ 1 - \left( P(W_{ij} < D_i)^m \times \int_0^{H_i} P\left(X_i(t) + y + u_i < H_i\right) f_{Y_i}^{<m>}(y) dy \right) \right] \times \frac{(\lambda_0 t)^m e^{(-\lambda_0 t)}}{m!} \quad (6)$$

## 4. Maintenance modeling

In this paper, we propose dynamically changing inspection planning as the preventive maintenance for series and parallel systems. It is assumed that the whole system and all the components are inspected at any inspection time, and each failed component is replaced at the beginning of the next inspection individually, while all other components continue functioning. Following that, the next inspection interval should be calculated based on the initial age of all the components. A penalty cost is assumed as the production loss for system downtime when the whole system is down. To find the optimal inspection interval for the next inspection Cost rate function should be calculated and optimized dynamically.

$$CR(\tau;\mathbf{u}) = \frac{C_I + C_R(1-R(\tau;\mathbf{u})) + C_\rho E[\rho]}{\tau} \quad (7)$$

Where $C_I$ is the inspection cost, $C_R$ is the cost of replacement, $C_\rho$ is the downtime cost and $\tau$ is the inspection time. $E[\rho]$ is the expected of downtime which is $E[\rho] = \int_0^\tau (\tau-t) f_{T_{H-U}}(t)dt$, where $f_{T_{H-u}}(t)$ is the probability density function of failure time that system fails at time $t$ during the time interval $\tau$ starting from random value as initial degradation $u$. it is calculated by taking derivative of CDF of system failure time $f_{T_{H-u}}(t) = \frac{dF_{T_{H-u}}(t)}{dt} = \frac{d(1-R(t;\mathbf{u}))}{dt}$.

According to the system configuration, equation (5) for series and equation (6) for parallel, should be substituted in equation (7) to calculate the system cost rate function. To find the next optimal inspection interval, an optimization problem should be solved dynamically based on the initial age of all the components, where cost rate function in equation (7) is the objective function of the optimization problem. Therefore, the next optimal inspection interval is obtaining dynamically based on initial random values as age of all the components.

## 5. Numerical example

In this paper two conceptual examples are considered to demonstrate the proposed reliability and maintenance model. The first example is a series system with three different components degrading with different rates. The second example is a parallel system with two components. Each component experiences two competing failure processes, soft and hard failure. Table 1 shows the parameter assumptions for these examples. It is assumed that the inspection cost is $5, replacement cost is $10, and downtime cost is $80. The cost rate function is calculated, and the optimization problem solved several times based on different combination of initial ages for series and parallel systems.

Table 1. Parameter Values for example

| Parameter | description | component 1 | component 2 | component 3 |
|---|---|---|---|---|
| $H_i$ | The soft failure threshold | 20 mm | 30 mm | 35 mm |
| $D_i$ | The hard failure threshold | 7 | 5 | 6 |
| $\alpha_i$ | The shape parameter for gamma process | 3 | 2 | 1 |
| $\beta_i$ | The scale parameter for gamma process | 1 | 0.6 | 0.3 |
| $\lambda$ | The initial intensity of random shocks | | $2.5\times 10^{-3}$ | |
| $Y_{ij}$ | Shock damage | $Y_{ij} \sim N(\mu_{Yi},\sigma_{YI}^2)$ $\mu_{Yi}=2$, $\sigma_{Yi}=0.5$ | $Y_{ij} \sim N(\mu_{Yi},\sigma_{YI}^2)$ $\mu_{Yi}=2.5$ GPa, $\sigma_{Yi}=0.2$ GPa | $Y_{ij} \sim N(\mu_{Yi},\sigma_{YI}^2)$ $\mu_{Yi}=3$ GPa, $\sigma_{Yi}=0.1$ GPa |
| $W_{ij}$ | The shock magnitude | $W_{ij} \sim N(\mu_{Wi},\sigma_{Wi}^2)$ $\mu_{Wi}=1.5$ $\sigma_{Wi}=0.4$ | $W_{ij} \sim N(\mu_{Wi},\sigma_{Wi}^2)$ $\mu_{Wi}=2$ GPa, $\sigma_{Wi}=0.3$ GPa | $W_{ij} \sim N(\mu_{Wi},\sigma_{Wi}^2)$ $\mu_{Wi}=1.2$ GPa, $\sigma_{Wi}=0.15$ GPa |

Figure 2 shows the cost rate function of three different combinations of components' initial age for the series example.

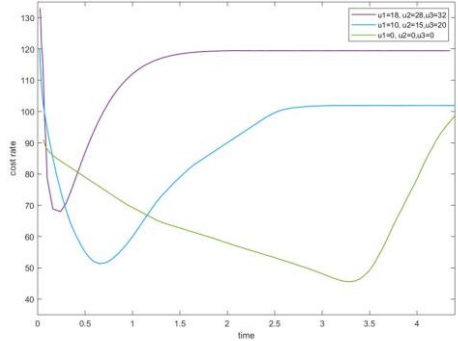

Figure 2. Cost rate for different combination of components' initial age

Table 2 and 3 show the optimal inspection interval for different combinations of initial ages for the series and parallel system.

As it is shown in Table 2. When all the three components are close to their failure threshold, the optimal inspection is very low, while when the initial ages are zero for all the components the optimal inspection is found as $\tau^* =3.3$. Moreover, it can be concluded that since component one degrades faster than other components, it has the higher impact on optimal inspection interval.

Comparing scenario 12 and 13, it can be concluded, that since the variance of component 3 is higher than component 2, it is dominant on determining optimal inspection interval. Scenarios 18-21 show than when one component or all the three components are close to their failure threshold the optimal inspection should be very low to replace them before failure.

Table 2. Optimal inspection interval for series system

| Scenarios number | Component 1 ($u_1$) | Component 2 ($u_2$) | Component 3 ($u_3$) | Optimal inspection interval ($\tau^*$) |
|---|---|---|---|---|
| 1 | 0 | 0 | 0 | 3.3 |
| 2 | 5 | 0 | 0 | 2.61 |
| 3 | 0 | 0 | 5 | 3.09 |
| 4 | 0 | 5 | 0 | 3.15 |
| 5 | 5 | 5 | 0 | 2.41 |
| 6 | 0 | 5 | 5 | 2.97 |
| 7 | 5 | 0 | 5 | 2.3 |
| 8 | 5 | 5 | 5 | 2.06 |
| 9 | 10 | 5 | 5 | 1.56 |
| 10 | 5 | 5 | 10 | 1.87 |
| 11 | 5 | 10 | 5 | 1.94 |
| 12 | 0 | 15 | 25 | 1.45 |
| 13 | 0 | 20 | 20 | 1.58 |
| 14 | 10 | 0 | 20 | 1.34 |
| 15 | 10 | 15 | 0 | 1.41 |
| 16 | 10 | 10 | 10 | 1.04 |
| 17 | 10 | 15 | 20 | 0.72 |
| 18 | 0 | 0 | 32 | 0.14 |
| 19 | 0 | 28 | 0 | 0.15 |
| 20 | 18 | 0 | 0 | 0.14 |
| 21 | 18 | 28 | 32 | 0.14 |

Table 3. Optimal inspection interval for parallel system

| Scenarios number | Component 1 ($u_1$) | Component 2 ($u_2$) | Optimal inspection interval ($\tau^*$) |
|---|---|---|---|
| 1 | 0 | 0 | 5.23 |
| 2 | 5 | 0 | 5.12 |
| 3 | 0 | 5 | 4.05 |
| 4 | 5 | 5 | 3.78 |
| 5 | 0 | 10 | 3.52 |
| 6 | 10 | 0 | 4.97 |
| 7 | 10 | 10 | 2.84 |
| 8 | 15 | 20 | 1.52 |
| 9 | 18 | 0 | 4.63 |
| 10 | 0 | 28 | 3.21 |
| 11 | 10 | 28 | 1.25 |
| 12 | 18 | 28 | 0.15 |

Table 3 shows the different optimal inspection interval for parallel example. Figure 3 also shows the result of Table 3 in a 3D plot. It is obvious that when the initial age of one of the components is zero, the optimal inspection is long. Moreover, when both components are close to their failure threshold, the optimal inspection is very short.

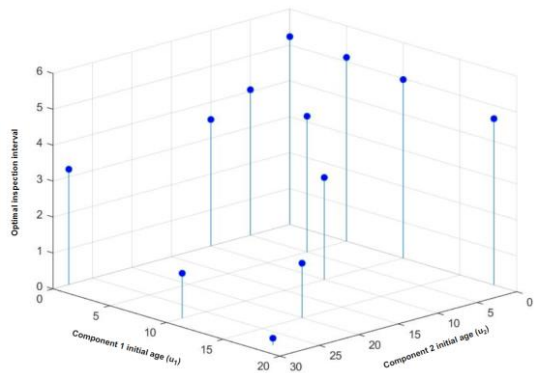

Figure 3. Combination of random age and inspection interval for parallel system

## 6. Conclusion

In this paper, a new maintenance model is developed for a system of multi-components where each component can be replaced individually within the system. It is considered that each component is subject to two dependent competing failure processes of soft and hard failure. gamma process models are used to model the degradation path of each component. Shock magnitude and shock damage are assumed to follow normal distributions. Reliability models are developed for series and parallel systems considering random or unspecified ages for each component at the beginning of an interval. The cost rate function is used to evaluate and optimize the maintenance model. It is conducted that the optimal inspection depends on the components' age at the beginning of that interval, so it is proposed that inspection interval should be found dynamically based on all the components' initial ages. Two numerical examples demonstrate how optimal inspection interval depend on the combination of random ages for series and parallel systems.

**References**

[1] Song, S., Coit, D. W., Feng, Q., & Peng, H. (2014). "Reliability analysis for multi-component systems subject to multiple dependent competing failure processes". *IEEE Transactions on Reliability*, 63(1), 331-345.
[2] Wang, Y., and H. Pham. "Dependent competing risk model with multiple-degradation and random shock using time-varying copulas." *16th ISSAT international conference on reliability and quality in design, Washington DC*. 2010
[3] Rafiee, K., Feng, Q., & Coit, D. W. (2014). "Reliability modeling for dependent competing failure processes with changing degradation rate". *IIE transactions*, 46(5), 483-496.
[4] Zhu, Y., Elsayed, E. A., Liao, H., & Chan, L. Y. (2010). "Availability optimization of systems subject to competing risk". *European Journal of Operational Research*, 202(3), 781-788.
[5] Song, S., Coit, D. W., & Feng, Q. (2014). "Reliability for systems of degrading components with distinct component shock sets". *Reliability Engineering & System Safety*, 132, 115-124.
[6] Abdul-Malak, David T., and Jeffrey P. Kharoufeh. "Optimally Replacing Multiple Systems in a Shared Environment." *Probability in the Engineering and Informational Sciences* 32.2 (2018): 179-206.
[7] Song, S., Yousefi, N., Coit, D. W., & Feng, Q. (2019). "Optimization of On-condition Thresholds for a System of Degrading Components with Competing Dependent Failure Processes". *arXiv preprint arXiv*:1901.07491.